\documentclass[12pt,draft]{article}
\usepackage{amssymb, amsmath, amsthm, verbatim, cite, lscape}
\usepackage[cp1251]{inputenc}
\usepackage[english]{babel}

\newtheorem{theorem}{Theorem}

\newcounter{parag}
\newcommand{\sect}[1]
{\refstepcounter{parag}
\begin{center} { \bf\theparag. #1} \end{center}}

\newtheorem{lemma}{Lemma}[parag]
\newtheorem{prop}[lemma]{Proposition}
\theoremstyle{definition}
\newtheorem*{prf}{Proof}

\theoremstyle{definition}

\begin{document}
\setcounter{page}{1} \thispagestyle{empty}
UDC 512.542

\begin{center}
\textbf{Locally finite groups of finite $c$-dimension\\~\\}
A.\,A. Buturlakin\end{center}

\begin{center}\textbf{Abstract}
\end{center}

{\small The supremum of lengths of strict chains of nested centralizers is called the $c$-dimension (centralizer dimension) of $G$. We prove two structure theorems for locally finite groups of finite $c$-dimension. We also prove that the $c$-dimension of the quotient of such a group by a locally soluble radical is bounded in terms of the $c$-dimension of the group itself.}

\begin{center}\textbf{Introduction}
\end{center}

Let $G$ be a group satisfying the minimal condition on centralizers ($\mathfrak M_c$-group). Following \cite{MS}, the supremum of lengths of strict chains of nested centralizers is called the $c$-dimension of $G$ and denoted by $\operatorname{cdim}(G)$. The main subject of this paper is locally finite groups of finite $c$-dimension. The class of groups of finite $c$-dimension includes some important classes of groups: linear groups, free groups, abelian groups. This class is closed under finite direct products, finite extensions, formation of subgroups, universal equivalence, but not closed under taking homomorphic images and extensions.

The class of locally finite $\mathfrak M_c$-groups possesses many interesting properties. For example, in~\cite{Bryant} it is proved that groups from this class satisfy Sylow's theorem for every prime and locally nilpotent $\mathfrak M_c$-groups are virtually nilpotent. In \cite{BryHar} locally soluble groups from this class was considered. The authors proved that these groups contain a nilpotent-by-abelian group of finite index and satisfy Hall's theorem on $\pi$-subgroups.

If one restricts his attention to groups of finite $c$-dimension, then it is natural to ask whether some quantative analogue of the results for $\mathfrak M_c$-groups can be proved. In \cite{Khukhro}, the author proved that the derived length of a locally soluble group of finite $c$-dimension is bounded in terms of $c$-dimension. He also gave some examples demonstrating that some direct analogues of corresponding result for $\mathfrak M_c$ do not hold: a soluble group of finite $c$-dimension $k$ does not possess a nilpotent-by-abelian subgroup of index bounded in terms of $k$, the rank of the factor of a nilpotent group of $c$-dimension $k$ by a subgroup of $k$-bounded nilpotency class is unbounded. In the same paper, a conjecture on the structure of arbitrary locally finite group of finite $c$-dimension was stated. Before formulating this conjecture, we should give some definitions.

The Hirsh--Plotkin radical $F(G)$ of a group $G$ is the maximal locally nilpotent subgroup of $G$. The component of $G$ is a subnormal quasisimple subgroup of $G$. The layer $E(G)$ is a subgroup generated by all components of $G$. The generalized Fitting subgroup $F^*(G)$ of a group $G$ is a subgroup generated by its Hirsh--Plotkin radical and layer.
\smallskip

\noindent\textbf{The Borovik-Khukhro conjecture.}
{\sl Let $G$ be a locally finite group of finite
$c$-dimension $k$. Let $S$ be the full inverse image of the generalized Fitting subgroup
$F^*(G/F(G))$. Then

$(1)$ the number of nonabelian simple composition factors of $G$ is finite and
$k$-bounded;

$(2)$ $G/S$ has an abelian subgroup of finite $k$-bounded index.}

\smallskip

While the first part of the was confirmed in \cite{13ButVas}, the second part does not hold even in the class of finite groups: the counterexample was constructed in \cite{BRV}. In the latter, a weaker version of the second part of the conjecture was proved for finite groups: it was proved that the conjecture becomes true if we replace $F(G)$ in the definition of $S$ be the soluble radical. In view of this result, A.V. Vasil'ev formulated the following question \cite[Question 19.101]{Ma}.

\smallskip

\noindent\textbf{Question. }\emph {Let $G$ be a locally finite group of finite $c$-dimension $k$, and let $S$ be
the preimage in $G$ of the socle of $G/R$, where $R$ is the locally soluble radical of $G$.
Is it true that the factor group $G/S$ contains an abelian subgroup of index bounded
by a function of $k$?}

\smallskip

We answer this question in the affirmative.

\begin{theorem}\label{t:Structure} Let $G$ be a locally finite group of $c$-dimension $k$ and $R$ be its (locally) soluble radical. Denote by $\overline G$ the quotient of $G$ by $R$ and by $H$ the socle of $\overline G$.  Then $C_{\overline G}(H)=1$ and the following statements hold.

$R$ has $k$-bounded derived length;

$H$ is a direct product of linear simple groups with $\lambda(H)$ bounded by a linear function of $k$;

$\overline G/H$ contains a finite abelian group $A$ with $l(A)$ and index bounded by a function of $k$.
\end{theorem}

Here $l(G)$ denotes the maximal length of strict chain of nested subgroup of a group $G$. The precise definition of $\lambda(G)$ is given below in Section 1  (without going into details, $\lambda(G)$ is the sum of Lie ranks of directs factors of~$S$).

We should note that a preprint with independent proof of similar result appeared on Arxiv recently \cite{Borovik}.

One of the tools in the proof of this theorem is the following statement.

\begin{theorem}\label{t:main} Let $G$ be a locally finite group of finite $c$-dimension $k$. Let $\overline G$ be the quotient of $G$ by its locally soluble radical. The $c$-dimension of $\overline G$ is bounded by a function of $k$.
\end{theorem}

We also prove another assertion about the structure of locally finite group of finite $c$-dimension in the final part of the paper, whose formulation is more similar to the one of the original Borovik--Khukhro conjecture (see, Theorem~\ref{t:Fitting}).

\sect{Preliminary results}

A group is called an $\mathfrak{M}_c$-group, if it satisfies the minimal condition on centralizers. The $c$-dimension $\operatorname{cdim}(G)$ of a group $G$ is the supremum of lengths of chains of nested centralizers. Denote by $l(G)$ the maximal length of a chain of nested subgroups of a group $G$.

It is known that a group $G$ possessing an $\mathfrak{M}_c$-subgroup $H$ of finite index is an $\mathfrak{M}_c$-group itself (the proof of this fact is similar to the proof of \cite[Lemma 3.20]{KegWeh}. In \cite[Proposition 3.8]{DKR}, it was proved that if additionally $c$-dimension of $H$ is finite and equal to $d$ and $|G:H|$ is equal to $k$, then $c$-dimension of $G$ is at most $k(k(d+2)+2)$. The proof of the following statement is similar to the proof of \cite[Lemma 3.20]{KegWeh}.

\begin{prop}\label{p:FinExt} Let $G$ be a group and $N$ be a normal subgroup of $G$. Assume that $l(G/N)$ is finite and is equal to $l$. If $N$ is an $\mathfrak{M}_c$-group, then $G$ is also an $\mathfrak{M}_c$-group. If $N$ has finite $c$-dimension $k$, then the $c$-dimension of $G$ is also finite and is at most $(l+1)^2(k+1)$.
\end{prop}

\begin{prf}  Let $\{C_i\}_{i\geqslant 1}$ be a properly ascending chain of centralizers. Since $l(G/N)$ is finite, the series $\{C_iN\}$ contains at most $l+1$ distinct subgroups. Therefore we may assume that $$C_iN=C_jN$$ for all subgroups of the series. Hence $$C_i=C_1(C_i\cap N)$$ and  $$C_G(C_i)=C_G(C_1)\cap C_G(C_i\cap N).$$ It follows that $\{C_G(C_i\cap N)\}_{i\geqslant1}$ is a strictly descending chain of centralizers. Since every centralizer $C_G(C_i\cap N)$ is an extension of some subgroup of $G/N$ by $C_N(C_i\cap N)$, the length of this series as well as of original series is finite. This proves the assertion for the case of $\mathfrak{M}_c$-groups. The same arguments show that if $\operatorname{cdim}(N)=k$, then $\operatorname{cdim}(G)\leqslant (l+1)^2(k+1)$.
\end{prf}

For a finite group $G$, we define a non-negative integer $\lambda(G)$ as follows. First, let $G$ be a nonabelian simple group. If $G$ is a group of Lie type, then $\lambda(G)$ is the minimum of Lie ranks of groups of Lie type isomorphic to $G$ (the Lie rank is the rank of the corresponding $(B,N)$-pair \cite[page 249]{Asch}). If $G$ is an alternating group, then let $\lambda(G)$ be the degree of $G$ (except for the groups $\operatorname{Alt}_5$, $\operatorname{Alt}_6$ and $\operatorname{Alt}_8$ which are isomorphic to groups of Lie type and, therefore, already have assigned values of $\lambda(G)$). Put $\lambda(G)=1$ for the sporadic groups. Now if $G$ is an arbitrary finite group, then $\lambda(G)$ is the sum of $\lambda(S)$ where $S$ runs over the nonabelian composition factors of $G$.

\begin{lemma}\label{l:Quasisimple}\emph{\cite[Proposition 3.1]{BRV}} There exists a universal constant $b$ such that ${\lambda(E(G))\leqslant b\cdot \operatorname{cdim}(G)}$ for every finite group $G$.
\end{lemma}

\begin{prop}\label{p:CommonRank}\emph{\cite[Proposition 3.2]{BRV}} There exists a universal constant $d$ such that $\lambda(G)\leqslant d\cdot k$ for every finite group $G$.
\end{prop}

\begin{lemma}\label{l:cdimGL} Let $G$ be a general linear group of dimension $n$ over a field. Then $\operatorname{cdim}(G)\leqslant n^2+1$.
\end{lemma}

\begin{prf} See the proof of \cite[Proposition 2.1]{MS}.
\end{prf}

Since a group of permutation matrices is isomorphic to corresponding symmetric group, the following lemma is an immediate corollary of the previous one.

\begin{lemma}\label{l:cdimAlt} $\operatorname{cdim}(A_n)\leqslant n^2+1$.
\end{lemma}

\begin{lemma}\label{l:cdimLie} Let $G$ be a group of Lie type of Lie rank $n$. The $c$-dimension of $G$ is bounded by a function of $n$.
\end{lemma}

\begin{prf} A group of Lie type is a group of automorphisms of some finite dimensional Lie algebra. The dimension of this algebra is a function of $n$. Therefore any group of Lie type can be embedded into $GL_{l}(q)$ for some $l$ depending on $n$ and the Lie type of $G$.  Since the $c$-dimension of a subgroup cannot exceed the $c$-dimension of the group, the lemma follows from Lemma~\ref{l:cdimGL}.
\end{prf}

\begin{lemma}\label{l:Khukhro}{\rm \cite[Lemma 3]{Khukhro}} If an elementary abelian $p$-group $E$ of order~$p^n$ acts faithfully on a finite nilpotent $p'$-group $Q$, then there exists a series of subgroups $$E=E_0>E_1>E_2>\dots>E_n=1$$ such that $$C_Q(E_0)<C_Q(E_1)<\dots<C_Q(E_n).$$
\end{lemma}

Let $H$ be a subgroup of $G$. The group $Aut_G(H)=N_G(H)/C_G(H)$ is called the group induced automorphisms of $H$ in $G$. Denote by $\overline G$ the quotient of $G$ by its soluble radical. Let the socle of $\overline G$ be a direct product of simple subgroups $S_1$, $\dots$, $S_k$. It is well known that $\overline G$ is isomorphic to a subgroup of a semidirect product of $Aut_{\overline{G}}(S_1)\times\dots\times Aut_{\overline{G}}(S_k)$ and some subgroup of the symmetric group $Sym_k$.

\begin{lemma}\label{l:IndAut} Let $G$ be a locally finite group ant $Q$ be its component. Let $R$ be the locally soluble radical of $G$. Denote by $\bar{\phantom{g}}$ the natural homomorphism from $G$ to $G/R$. Then $Aut_G(Q)\simeq Aut_{\overline G}(\overline Q)$.
\end{lemma}

\begin{prf} By the definition, $Aut_{\overline G}(\overline Q)\simeq N_{\overline G}(\overline Q)/C_{\overline G}(\overline Q)$. The full inverse image of $N_{\overline G}(\overline Q)$ in $G$ is equal to $N_G(QR)R$. From the equality $Q=E(QR)$ it follows that $Q$ is characteristic in $QR$, and therefore $N_G(QR)R=N_G(Q)$. Since $\overline {C_G(Q)}\leqslant C_{\overline G}(\overline Q)$, the group $Aut_{\overline G}(\overline Q)$ is a homomorphic image of $Aut_G(Q)$. But both of them are almost simple, hence they are isomorphic as required.
\end{prf}

\begin{lemma}\label{l:FinInInfin} Let $G$ be a locally finite group of finite $c$-dimension $k$. Then $G$ contains a finite subgroup $H$ such that $cdim(H)=k$.
\end{lemma}

\begin{prf} Since $c$-dimension of $G$ is equal to $k$, there exist $x_1$, $\dots$, $x_k$ in $G$ such that $$C_G(x_1)>C_G(x_1, x_2)>\dots>C_G(x_1, x_2,\dots, x_k)$$ is a centralizer chain of maximal length. Chose $h_i\in C_G(x_1, x_2,\dots, x_i)\setminus C_G(x_1, x_2,\dots, x_{i+1})$ for $1\leqslant i<k$. It is easy to see that the $c$-dimension of a finite group generated by $x_1$, $\dots$, $x_k$ and $h_1$, $\dots$, $h_{k-1}$ is equal to $k$.
\end{prf}

For a group of Lie type $S$ over a field $k$, denote by $\Phi_S$ the group of field automorphisms of $S$, that is automorphisms acting on root subgroups by $x_r(t)\mapsto x_r(t^\sigma)$, where $\sigma$ is an automorphism of $k$.

Let us recall some information about locally finite fields. Basic information about locally finite fields and their automorphism groups is contained, for example, in \cite{BraSch}.

A Steinitz number is a formal product $$\prod\limits_{i=0}^\infty p_i^{x_i},$$ where $p_i$ denotes the $i$th prime number and $x_i$ is an element of the extended set of natural numbers $\{ 0, 1, \dots, n, \dots, \infty\}$. There are natural embedding of positive integer into the set of Steinitz numbers and extension of divisibility relation to this set. Let $N$ be a Steinitz number and $q$ a prime power. Put $$GF(q^N)=\bigcup\limits_{d|N}GF(q^d),$$ where union is taken over all ordinary natural divisors $d$. Every locally finite field is isomorphic to $GF(q^N)$ for some Steinitz number $N$.

The automorphism group of $GF(q^N)$ is the inverse limit of automorphism groups of its finite subfields. Every finite homomorphic image of $Aut(GF(q^N))$ is cyclic.

\begin{lemma}\label{l:ExtByField} Let $G$ be an group with the layer isomorphic to a group of Lie type $S$ over a locally finite field $K$, such that $G/S\leqslant \Phi_S$. If the $c$-dimension of $G$ is finite and equal to $k$, then $l(G/S)\leqslant k$.
\end{lemma}

\begin{prf} Assume that $l(G/S)>k$. Hence there is a finite subfield $k$ of $K$, such that the image $A$ of $G/S$ in the automorphism group of $k$ satisfies $l(A)>k$. Let $X$ be a subgroup of a root subgroup consisting of elements $x_r(t)$ where $t\in k$. If $B$ and $C$ are distinct subgroups of $A$, then $C_X(B)\neq C_X(C)$. Therefore, the $c$-dimension of $G$ is at least $l(A)$, that is a contradiction.
\end{prf}

In the following lemma $S_p(G)$ stands for the locally $p$-soluble radical of $G$. A group is called $p$-perfect, if it is perfect and generated by its $p$-elements. A group satisfies the strong Sylow theorem for a prime $p$, if all its subgroups satisfy the Sylow theorem for this prime.

\begin{lemma}\label{l:Kegel}\emph{\cite[Theorem 4.3]{Kegel}} If the locally finite group $G$ is not $p$-soluble and satisfies the strong Sylow theorem for the prime $p\geqslant 5$, then the socle $Soc(X)$ of the factor group $X\simeq G/S_p(G)$ is the direct product of finitely many linear simple $p$-perfect subgroups. Also the centralizer $C_X(Soc(X))$ is trivial.
\end{lemma}

\sect{Proof of Theorem \ref{t:main}}

We start by proving the theorem in the finite case.

\begin{prop}\label{p:main} Let $G$ be a finite group of finite $c$-dimension $k$. Let $\overline G$ be the quotient of $G$ by its soluble radical. The $c$-dimension of $\overline G$ is bounded by a function of $k$.
\end{prop}

\begin{prf}  The socle of $\overline G$ is a direct product of finite nonabelian simple groups. From Proposition~\ref{p:CommonRank}, it follows that $\lambda(Soc(\overline G))$ is bounded by a linear function of $k$; in particular, the number of composition factors of $Soc(\overline G)$ is $k$-bounded. By Lemmas~\ref{l:cdimAlt} and \ref{l:cdimLie}, the $c$-dimension of $Soc(\overline G)$ is in terms of $k$. By Proposition~\ref{p:FinExt} the $c$-dimension of $\overline G$ is bounded in terms of the $c$-dimension of $Soc(\overline G)$ and $l(\overline G/Soc(\overline G))$. Hence it is sufficient to bound the latter by a function of $k$. Let us show that $l(Aut_{\overline G} S)$ is $k$-bounded for every composition factor of $Soc(\overline G)$. Since the number of composition factors of $Soc(\overline G)$ is $k$-bounded, this will prove the proposition.

Put $\Phi=\prod_S \Phi_S$, where $S$ ranges through all composition factors of $Soc(\overline G)$. By the classification of finite simple groups, the index of $\Phi$ in $Out(Soc(\overline G))$ is bounded in terms of $\lambda(Soc(\overline G))$. Hence the index of $\Phi\cap \overline G/Soc(\overline G)$ in $\overline G/Soc(\overline G)$ is bounded in terms of $k$. Therefore, by Proposition~\ref{p:FinExt}, we may assume that the quotient $\overline G/Soc(\overline G)$ is a subgroup of $\Phi$.

Let $R$ be the soluble radical of $G$. By Frattini's argument, we may assume that $R$ is nilpotent. Let $S$ be a composition factor of $Soc(\overline G)$ which is a group of Lie type over a field of characteristic $p$ and order $p^\alpha$. If $X_r$ is a root subgroup of $S$, then the center $Z(X_r)$ is an elementary abelian group whose rank is at least $\alpha/2$. The group $Z(X_r)$ acts on $O_{p'}(R)$, and either this action is faithful, or $S$ is a composition factor of $C_G(O_{p'}(R))$. If the action is faithful, then Lemma~\ref{l:Khukhro} yields $\alpha/2\leqslant k$. Therefore the order of $\Phi_S$ is at most $2k$.

Assume that $S$ is a composition factor of $C_G(O_{p'}(R))$. We again have two cases: either $S$ is a composition factor of $C_G(R)$, or $S$ is a composition factor of $Aut_G(O_p(R))$. If $S$ is a composition factor of $C_G(R)$, then $S$ is a composition factor of the layer $E(G)$ and there is a quasisimple subgroup $T$ of $G$ such that $T/Z(T)\simeq S$. By Lemma~\ref{l:IndAut}, the group $Aut_G(T)$ is isomorphic to $Aut_{\overline G}(S)$. By Lemma~\ref{l:ExtByField} we have $l(Aut_G(T))\leqslant k$ as required.

Let $S$ be a composition factor of $Aut_G(O_p(R))$. Since $S$ is a composition factor of $C_G(R)$, it is a composition factor of the layer of $G/O_p(R)$. Denote by $T$ corresponding component of $G/O_p(R)$. With the finite number of exceptions of isomorphism types of $S$ (for the exact list of exceptions see, for example, \cite[Table 6.1.3]{GLS}) which can be omitted from the considerations, the group $T$ is also a group of Lie type. By Lemma~\ref{l:IndAut}, the group of induced automorphisms of $T$ coincides with $Aut_{\overline G}(S)$. By Lemma~\ref{l:ExtByField} we have $l(Aut_{\overline G}(S)/S)$ is at most $k$

Put $A=Aut_{\overline G}(S)/S$ ($A$ can be assumed to be a subgroup of $\Phi_S$ ). Consider a Cartan subgroup $C$ of $T$ (its structure is described in \cite[Theorem 2.4.7]{GLS}). If $\sigma_1$ and $\sigma_2$ are two automorphisms from $\Phi_S$ of distinct orders, which are automorphisms of $T$, then $C_{C}(\sigma_1)\neq C_{C}(\sigma_2)$. Hence for each cyclic subgroup of $A$ there is an element of $C$ centralized by this subgroup, but not larger. Therefore there are elements $c_1$, $c_2$, $\dots$, $c_{l(A)}$ of $C$ such that $$C_{G/O_p(R)}(c_1, c_2,\dots, c_i)>C_{G/O_p(R)}(c_1, c_2,\dots, c_{i+1})$$ for $1\leqslant i<l(A)$. By the properties of coprime action, there exist inverse images $c_1'$, $\dots$, $c_{l(A)}'$ of these elements in $G$ such that $$C_{G}(c_1, c_2,\dots, c_i)>C_{G}(c_1, c_2,\dots, c_{i+1})$$ for $1\leqslant i<l(A)$. Hence $l(A)\leqslant cdim(G)$, as was mentioned before, this completes the proof of the proposition.
\end{prf}

\noindent\textbf{Proof of Theorem \ref{t:main}.} Let $f$ be a function such that, for any finite group $G$ and its soluble radical $R$, we have $$cdim(G/R)\leqslant f(cdim(G)).$$ We will assume that $f$ is nondecreasing. Let $G$ be a locally finite group of $c$-dimension $k$ and $R$ be its locally soluble radical. Assume that $$cdim(G/R)>f(cdim(G)).$$ By Lemma~\ref{l:FinInInfin}, there is a finite group $\overline H$ in $G/R$ such that $cdim(H)>f(cdim(G))$. Let $H$ be a a finite preimage of $\overline H$ in $G$. Let $\Sigma$ be a local system of $G$ consisting of finite groups (for the definition, see \cite[page 8]{KegWeh}). There exists an element $X$ of $\Sigma$ containing $H$. It is known (see \cite[pp. 12-15]{KegWeh}) that there is $Y\in\Sigma$ such that $X\cap R$ is equal to the intersection of $X$ and the soluble radical $R(Y)$ of $Y$. Hence $\overline H$ is isomorphic to a subgroup of $Y/R(Y)$. Since $$cdim(Y)\leqslant cdim(G),$$ this a contradiction.

\sect{Proof of Theorem \ref{t:Structure}}

The assertion about the locally soluble radical is just a part of the main theorem of \cite{Khukhro}. The proof of the statements about the socle of $G/R$ is collected in the following proposition.

\begin{prop}\label{p:Kegel} Let $G$ be a locally finite group of $c$-dimension $k$ and $R$ be its locally soluble radical. Then the socle of $\overline G=G/R$ is the direct product of finitely many linear simple groups the sum of Lie ranks of which is bounded by a linear function of $k$. Also the centralizer $C_{\overline G}(Soc(\overline G))$ is trivial.
\end{prop}

\begin{prf} By the main theorem of \cite{13ButVas}, the number of nonabelian composition factors of $G$ is less than $5k$. Hence there is a finite set of primes $P$ such that every prime in $P$ is greater than $3$ and every nonabelian composition factor of $G$ contains an element whose order lies in $P$.

From \cite[Theorem B]{Bryant} it follows that a locally finite group with minimal condition on centralizers satisfies the strong Sylow theorem for every prime~$p$. Therefore due to Lemma~\ref{l:Kegel}, the group $\overline G$ is isomorphic to a subdirect product of the groups $X_p=G/S_p(G)$ for $p\in P$. Moreover, $C_{X_p}(Soc(X_p))=1$ for every $p\in P$.

First, let us prove that $C_{\overline G}(Soc(\overline G))=1$. Since the locally soluble radical of $\overline G$ is trivial, it is sufficient to prove that every nontrivial normal subgroup $N$ of $\overline G$ contains a minimal normal subgroup. Denote by $\varphi_p$ the projection map on $X_p$. We can assume that $N\varphi_p$ is either trivial, or a minimal normal subgroup of $X_p$. Indeed, if $N\varphi_p\neq 1$, then it contains a minimal normal $M_p$ subgroup of $X_p$ and, without loss of generality, we may replace $N$ by the full inverse image of $M_p$ in $N$. Therefore, $N$ is a sibdirect product of finite number of simple groups and is itself a direct product of simple groups. Hence $N$ contains a minimal normal subgroup of $G$. By using similar arguments, one can show that every minimal subgroup of $\overline G$ is a direct product of finite number of simple groups, so the same is true for $Soc(\overline G)$.

Every linear locally finite simple group is a group of Lie type over a locally finite field, and therefore contains a finite subgroup which is a group of the same Lie type. Hence the assertion about boundedness of the sum of Lie rank follows from Proposition~\ref{p:CommonRank}. This completes the proof of the proposition.
\end{prf}

\smallskip

\textbf{Proof of Theorem \ref{t:Structure}}

By Theorem~\ref{t:main}, we may assume that $R=1$. Hence, $H$ is equal to $Soc(G)$. Since $C_G(H)=1$, we may regard $G$ as a subgroup of $Aut(H)$. Let $\Phi$ be the product of $\Phi_S$ where $S$ runs through all composition factors of $H$. Let $F$ be the full inverse image of $\Phi$ in $Aut(H)$. Put $K=G\cap F$. Since the index of $F$ in $Aut(H)$ is bounded in terms of $\lambda(H)$, the index of $K$ in $G$ is bounded. Also $K/H$ is abelian. Hence it remains to show that $l(K/H)$ is $k$-bounded. From Lemma~\ref{l:ExtByField} it follows that $l(Aut_K(S)/S)$ is bounded by a function of $k$ for every composition factor $S$ of $H$. Since the number of composition factors of $H$ is also bounded in terms of $k$, this completes the proof of the theorem.

As a corollary of Theorem~\ref{t:Structure}, we obtain the following assertion.

\begin{theorem}\label{t:Fitting} Let $G$ be a locally finite group of $c$-dimension $k$. Let $\overline G$ be the quotient group of $G$ by the third Hirsch--Plotkin radical $F_3(G)$. Then the quotient of $\overline G$ by the layer $E(\overline G)$ contains a finite abelian group $A$ with $l(A)$ and index bounded by a function of $k$.
\end{theorem}

\begin{prf} By the main theorem of \cite{Khukhro}, the locally soluble radical $R$ of $\overline G$ is finite of order bounded by a function of $k$.

Let us prove that the centralizer  $C_{\overline G}(F^*(\overline G))$ is contained in $F^*(\overline G)$. Put $C=C_{\overline G}(F^*(G))$. It is easy to see that $F^*(C)\leqslant F^*(\overline G)\leqslant C_{\overline G}(C)$. Therefore, $F^*(C)=Z(C)$. Assume that $C\neq Z(C)$. If $C\cap R\neq Z(C)$, then $C/Z(C)$ contains a cyclic subnormal subgroup whose full inverse image is abelian and subnormal in $C$; that is a contradiction. Hence $C\cap R=Z(C)$. Therefore, $C/Z(C)$ is a nontrivial normal subgroup of $\overline G/R$. Theorem~\ref{t:Structure} implies that the socle of $\overline G/R$ is the direct product of finitely many simple groups and the centralizer of the socle is trivial. This means that $C\cap Soc(\overline G/R)\neq 1$ and $C$ contains a minimal subnormal subgroup with is a nonabelian simple. Let $E$ be its full inverse image in $C$. Then we have $E=E'Z(C)$ where $E$ is a component of $C$; that is a final contradiction.

Since $C\leqslant F^*(\overline G)$, the group $\overline G/C$ is contained in the direct product of $Aut_{\overline G}(E(\overline G))$ and $Aut_{\overline G}(F(\overline G))$. Since the order of $F(\overline G)$ is bounded in terms of $k$, the order of $Aut_{\overline G}(F(\overline G))$ is also bounded. Therefore, we may consider only $Aut_{\overline G}(E(\overline G))$. If $C_1$, $\dots$, $C_s$ are components of $\overline G$, then $Aut_{\overline G}(E(\overline G))$ is a subgroup of semidirect product of $$Aut(C_1)\times\dots Aut(C_s) \text{ and a subgroup of } Sym_s.$$ Since $s$ is bounded in terms of $k$, we may assume that $Aut_{\overline G}(E(\overline G))$ is a subgroup of direct product of automorphism groups of components, that is $Aut_{\overline G}(E(\overline G))$ is a subgroup of the direct product of $Aut_{\overline G}(C_i)$ for $1\leqslant i\leqslant s$. Now the theorem follows from Lemma~\ref{l:IndAut} and Theorem~\ref{t:Structure}.
\end{prf}

\end{document}